\documentclass{article}

\usepackage{amsfonts}
\usepackage{amsmath}
\usepackage{graphicx}
\usepackage{psfrag}
\numberwithin{equation}{section}                                

\newtheorem{prop}{Proposition}
\newtheorem{thm}{Theorem}
\newtheorem{rem}{Remark}

\newcommand{\proof}{\smallskip\par\noindent\emph{Proof.}\quad}
\newcommand{\qed}{\nobreak\quad\nobreak\hfill\nobreak$\square$\vspace{8pt}\par}

\newcommand{\R}{\mathbb{R}}
\newcommand{\Z}{\mathbb{Z}}
\newcommand{\n}{\noindent}
\renewcommand{\eqref}[1]{(\ref{#1})}
\newcommand{\disp}{\displaystyle}
\newcommand{\Og}{{\mathcal O}}
\newcommand{\I}{{\mathcal I}}
\newcommand{\ini}{{\rm in}}

\title{Infinite Horizon Noncooperative Differential Games with Non-Smooth Costs}
\author{Fabio S. Priuli}
\date{}

\begin{document}
\maketitle

\section{Introduction}

This paper deals with the study of a class of non-cooperative
differential games in infinite time horizon. Namely, we consider a
game with dynamics

\begin{equation}
\dot x=\sum_{i=1}^n\alpha_i\,,\qquad\qquad x(0)=y\,,
\end{equation}

\n where each player acts on his control $\alpha_i$ to minimize an
exponentially discounted cost of the form

\begin{equation}\label{eq:costs_intro}
J_i(\alpha)\doteq\int_0^\infty e^{-t} \bigg[ h_i\big(x(t)) +
{\alpha_i^2(t)\over 2} \bigg]\,dt\,,
\end{equation}

\n both $h_i$ being integrable functions, whose smoothness will be
addressed later. Very few results are known on the subject, except
in two particular cases: two players zero-sum games and LQ games
(where LQ stands for linear-quadratic). Indeed, a key step in this
kind of problems is the study of the value function $u$. In the
region where $u$ is smooth, its components satisfy a system of
Hamilton-Jacobi equations (see~\cite{FriedmanBook}), and this system
is usually difficult to solve.

In the case of two players zero-sum games, since what one player
gains is exactly what the other player loses, the two components of
$u$ are one the opposite of the other. Hence, the Hamilton-Jacobi
system (HJ in the following) reduces to a single equation and one
can apply the standard theory of viscosity solutions
(see~\cite{BardiCapDolcetta} for more details) to obtain existence
and uniqueness results.

In the case of LQ games, the HJ system can be connected to a Riccati
system of ODE for matrices. This system is in general much easier
than the original one, and standard ODE techniques can be applied
(see~\cite{EngwerdaBook} for a detailed treatment).

On the other hand, in the present case both approaches fail and
therefore one has no established techniques to rely on. However, few
results still can be proved.

In the finite horizon setting, the analysis presented
in~\cite{BressanShen1,BressanShen2} showed that, for a
non-cooperative $n$-players differential game with general terminal
payoffs, the well-posedness is strongly related with the HJ system
being hyperbolic. Namely, for $x\in\R$, thanks to recent advances in
the theory of hyperbolic systems of PDE, games with strictly
hyperbolic HJ systems are well-posed. On the other hand, it is
possible to produce examples of games, even in one spatial
dimension, whose corresponding HJ system is not even weakly
hyperbolic and, hence, it is ill-posed.

A first attempt to study this problem in the infinite horizon
setting, for two players, was made in~\cite{BressanPriuli}. The same
simple game was considered, and it was proved that, depending on the
monotonicity of the cost functions, very different situations could
arise. Indeed, the HJ system in this case takes the following form

\begin{equation}\label{eq:HJ_first}
\left\{\begin{array}{l}
u_1(x)=h_1(x)-u_1'u_2'-(u_1')^2/2\,,\\
\\
u_2(x)=h_2(x)-u_1'u_2'-(u_2')^2/2\,.\\
\end{array}\right.
\end{equation}

But with a system of this form, we can end up with too many
solutions. We find not only value functions $u$ that leads to Nash
equilibria in feedback form, but also solutions that does not
represent equilibria of the game. It is then necessary to introduce
a suitable concept of admissibility. In particular we say that a
solution $u$ is {\it admissible}, if $u$ is a Carath\'eodory
solution of \eqref{eq:HJ_first}, which grows at most linearly as
$|x|\to\infty$ and satisfies suitable jump conditions in points
where its derivatives are discontinuous. For such a kind of
solutions, a {\it verification theorem} was proved: given an
admissible solution $u$ and denoted by $u'_i$ the components of its
derivatives, then $\alpha_i=-u'_i$ provide a Nash equilibrium
solution in feedback form.

In~\cite{BressanPriuli}, it turned out that existence and uniqueness
of admissible solution for \eqref{eq:HJ_first} heavily depend on the
choice of the costs.

First, suppose that both the cost functionals are increasing (resp.
decreasing). This means that both players would like to steer the
game in the same direction, namely the direction along which their
costs decreases. In this case an admissible solution always exists,
and it is also unique, provided a {\em small oscillations}
assumption is satisfied. This existence result was in some sense
expected, since this case corresponds, in the finite horizon
setting, to the hyperbolic one studied in~\cite{BressanShen1}.

Suppose now that the cost functionals have opposite monotonicity.
This means that the players have conflicting interests, since they
would like the game to go in different directions. In this case it
is known, see~\cite{BressanShen2}, that the finite horizon problem
is in general ill-posed. On the same line, for our game, it is
enough to consider two linear functionals with opposite slopes (say
$k, -k$, for any real number $k\neq 0$) to find infinitely many
admissible solutions, and hence infinitely many Nash equilibria in
feedback form. Nevertheless, quite surprisingly, it's still possible
to recover existence and uniqueness of admissible solutions to
\eqref{eq:HJ_first} in the case of costs that are small perturbation
of linear ones, but with slopes that are not exactly opposite.

This richness of different situations reflects in some sense the
results found in~\cite{CardPlas}. Indeed, the exact same dynamics
was studied, in the finite horizon case, with only exit costs. Main
differences between~\cite{CardPlas} and~\cite{BressanPriuli,
BressanShen1, BressanShen2} lay in the concept of solution. The
authors of~\cite{CardPlas} look for discontinuous feedback controls
that not only leads to Nash equilibria, but also satisfies a sort of
programming principle. This resulted in (uncountable) infinitely
many solutions, at price of stronger assumptions on the final costs.

While the cost functionals considered in \cite{BressanPriuli} were a
small perturbation of affine costs, in the present paper we study a
wider class of cost functions. Motivated by the theory of hyperbolic
systems~\cite{BressanBook}, we now consider piecewise linear cost
functionals, whose derivative has jumps. This setting is a natural
first step towards the analysis of existence and uniqueness of Nash
equilibrium solutions for non-linear costs.

Again, as in~\cite{BressanPriuli}, we reach different results
depending on the signs chosen for $h'_i$. Indeed, as it will be
proved in the following sections, if we are in the cooperative
situation for all $x$, we can still recover a unique admissible
solution for~\eqref{eq:HJ_first}. On the other hand, any change in
the behavior of the costs will translate in some sort of instability
of the game, leading either to infinitely many admissible solution,
or to one unique admissible solution, or even to no admissible
solution at all, only depending on the particular choices of the
slopes $h'_i$.

In conclusion, this great variety of arising situations seems to
suggest that the present approach is not the most suitable one to
deal with the intrinsic issues of the problem. In particular, we can
provide examples of very simple differential games where no
Carath\'eodory solution with sublinear growth at infinity exists.
Recalling that, in the case of smooth costs
(see~\cite{BressanPriuli}), this class of solutions was exactly the
right one to find Nash equilibria in feedback form, our study
strongly suggest that a different approach is needed: either to look
for Pareto optima, as in~\cite{BressanShen2}, or to introduce some
other relaxed concept of equilibrium.

The structure of the present paper is the following. In Section 2 we
will introduce main notations and definitions. Moreover we will
recall briefly what was proved in the case of smooth costs and
provide a couple of useful Lemmas. In Section 3 we will present and
prove the main results of this paper, dealing with cooperative
players, in the sense of players whose costs always have the same
monotonicity. In this case existence and uniqueness results hold for
both piecewise linear and piecewise smooth cost functionals. In
Section 4 we will prove that a similar extension is not possible in
the case of conflicting interests. Actually, we will provide an
example in which the games has infinitely many Nash equilibria, as
well as an example in which there cannot be any admissible solution
to~\eqref{eq:HJ_first}. Finally, in Section 5, we will discuss a
last case that can arise when either one or both the cost
functionals are allowed to change monotonicity. From this game which
is partially ``cooperative'' (in the sense above) and partially
``conflicting'', infinitely many Nash equilibria can be found.

\section{Preliminaries}

In this paper we consider a scalar $2$-persons differential game,
with dynamics

\begin{equation}\label{eq:dynamics}
\dot x=\alpha_1+\alpha_2\,,
\end{equation}
\begin{equation}\label{eq:init_state}
x(0)=y\,.
\end{equation}

\n The functions $t\mapsto \alpha_i(t)$, $i=1,2$, represent the
controls implemented by the $i$-th player, chosen within a compact
set of admissible controls ${\cal A}_i\subset\R$. The game takes
place on $[0,+\infty[$ and each player is subject to a running cost,
exponentially discounted, of the following form

\begin{equation}\label{eq:costs}
J_i(\alpha_i)\doteq\int_0^\infty e^{-t} \bigg[ h_i\big(x(t)) +
{\alpha_i^2(t)\over 2} \bigg]\,dt\,.
\end{equation}

\n Assume here that both $h_i$ are piecewise smooth functions with
bounded derivatives. Later we will weaken this requirements.

\n A couple of feedback strategies $(\alpha_1^*(x),\alpha_2^*(x))$
represents a {\em Nash equilibrium solution} for the
game~\eqref{eq:dynamics}--\eqref{eq:init_state} if the following
holds. For $i\in \{1,2\}$, the feedback control
$\alpha_i=\alpha_i^*(x)$ provides a solution to the the optimal
control problem for the $i$-th player,
\begin{equation}\label{eq:Nash1}
\min_{\alpha_i(\cdot)}~J_i(\alpha_i)\,,
\end{equation}
where the dynamics of the system is
\begin{equation}\label{eq:Nash2}
\dot x=f_i(x,\alpha_i)+f_j(x,\alpha^*_j(x)),
\qquad\qquad \alpha_i(t)\in A_i,\,j\neq i\,.\\
\end{equation}
More precisely, we require that, for every initial data $y\in\R$,
the Cauchy problem

\begin{equation}\label{eq:Nash3}
\dot x=
f_1\big(x,\alpha_1^*(x)\big)+f_2\big(x,\alpha_2^*(x)\big)\,,\qquad
\qquad x(0)=y\,,
\end{equation}

\n should have at least one Caratheodory solution $t\mapsto x(t)$,
defined for all $t\in [0,\infty[\,$. Moreover, for every such
solution and each $i=1,\ldots,m$, the cost to the $i$-th player
should provide the minimum for the optimal control
problem~\eqref{eq:Nash1}-\eqref{eq:Nash2}. We recall that a
Caratheodory solution is an absolutely continuous function $t\mapsto
x(t)$ which satisfies the differential equation in~\eqref{eq:Nash3}
at almost every $t> 0$.

By the theory of optimal control, see for
example~\cite{BardiCapDolcetta}, we know that if $u$ is the value
function corresponding to~\eqref{eq:dynamics}-\eqref{eq:init_state}
with costs
$$
J_i(\alpha_i)\doteq \int_0^\infty
e^{-t}\,\psi_i\big(x(t),\,\alpha_i(t)\big)\,dt\,,
$$
\n then, where $u$ is smooth, each component $u_i$ should provide a
solution to the corresponding scalar Hamilton-Jacobi-Bellman
equation. The vector function $u$ thus satisfies the stationary
system of equations
\begin{equation}\label{eq:HJ_general}
u_i(x)=H_i(x,\, u'_1, u'_2)\,,
\end{equation}
\n where the Hamiltonian functions $H_i$ are defined as follows. For
each $p\in\R$, assume that there exists an optimal control value
$\alpha_j^*(x,p)$ such that
\begin{equation}\label{eq:argmax}
p\cdot \alpha_j^*(x,p)+
\psi_j\big(x,\,\alpha_j^*(x,p)\big)=\min_{a\in A_j} \,\big\{p\cdot
a+ \psi_j(x,a)\big\}\,.
\end{equation}
Then
\begin{equation}\label{eq:hamiltonian}
H_i(x,\,p_1,p_2)\doteq
p_i\cdot\alpha_j^*(x,p_j)+\psi_i\big(x,\,\alpha_i^*(x,p_i)\big)\,.
\end{equation}

\n for $i,j\in\{1,2\}$ and $i\neq j$. In general, even in cases as
easy as $\psi_i=\alpha_i^2/2$, this system will have infinitely many
solutions defined on the whole $\R$ (see Example~1
in~\cite{BressanPriuli}). And not every solution corresponds to a
Nash equilibrium for the initial game. To single out a (hopefully
unique) admissible solution, and therefore a Nash equilibrium for
the differential game, additional requirements must be imposed.
Namely a solution $u$ to~\eqref{eq:HJ_general} is said to be an {\em
admissible solution} if the following holds:

\begin{description}
\item{{\bf (A1)}} ~$u$ is absolutely continuous and its derivative $u'$
satisfies~\eqref{eq:HJ_general} at a.e.~point $x\in\R$.

\item{{\bf (A2)}} ~$u$ has sublinear growth at infinity; namely,
there exists a constant $C$ such that, for all $x\in\R$,

\begin{equation}\label{eq:growth}
\big|u(x)\big|\leq C\,\big(1+|x|\big)\,.
\end{equation}

\item{{\bf (A3)}}  At every point $y\in \R$, the derivative $u'$ admits
right and left limits $u'(y+)$, $u'(y-)$ and at points where $u'$ is
discontinuous, these limits satisfy at least one of the conditions

\begin{equation}\label{eq:jumps}
{u_1'(y+)+u_2'(y+)}\leq 0\qquad\hbox{or}\qquad
{u_1'(y-)+u_2'(y-)}\geq 0\,.
\end{equation}
\end{description}

Because of the assumption on $h'_i$, the cost functions $h_i$ are
Lipschitz continuous. It is thus natural to require the value
functions $u_i$ to be absolutely continuous, with sub-linear growth
as $x\to\pm\infty$. The motivation for the assumption (A3) is quite
simple. Observing that, in~\eqref{eq:argmax}, the feedback controls
are $\alpha_i^*=-u_i'$, the condition~\eqref{eq:jumps} provides the
existence of a local solution to the Cauchy problem
$$
\dot x=-u_1'(x)-u'_2(x)\,,\qquad\qquad x(0)=y
$$
\n forward in time.  In the opposite case, solutions of the
O.D.E.~would approach $y$ from both sides, and be trapped.

Notice that, for $2$-players games, the assumptions (A3) is
equivalent to
\begin{equation}\label{eq:jumps_2p}
{u_1'(y+)}+{u_2'(y+)}\leq 0\,,\qquad \qquad u_i'(y-)= -u_i'(y+)
\qquad (i=1,2)\,.
\end{equation}

This concept of admissibility turns out to be the right one. Indeed,
the following verification theorem can be proved (see
again~\cite{BressanPriuli}).

\begin{thm}
Consider the differential
game~\eqref{eq:dynamics}--\eqref{eq:init_state}. Let $u:\R\mapsto
\R^m$ be an admissible solution to the systems of H-J
equations~\eqref{eq:HJ_general}, so that the conditions (A1)--(A3)
hold.  Then the controls $\alpha_i^*=-u'_i$ provide a Nash
equilibrium solution in feedback form.
\end{thm}

Anyway, this theorem says nothing about the actual existence of
admissible solutions to~\eqref{eq:HJ_general}. To deal with this
problem, some manipulations have to be done on~\eqref{eq:HJ_general}
itself. Indeed, in the present case of costs as
in~\eqref{eq:costs_intro}, the Hamiltonian
functions~\eqref{eq:hamiltonian} lead to

\begin{equation}\label{eq:HJ_impl}
\left\{\begin{array}{l}
u_1(x)=h_1(x)-u_1'u_2'-(u_1')^2/2\,,\\
\\
u_2(x)=h_2(x)-u_1'u_2'-(u_2')^2/2\,.\\
\end{array}\right.
\end{equation}

\n Differentiating~\eqref{eq:HJ_impl} w.r.t. $x$ and setting
$p_i=u_i'$ one obtains the system
\begin{equation}
\left\{\begin{aligned}
h_1'-p_1&= (p_1+p_2)p_1'+p_1p_2'\,,\\
&\\
h_2'-p_2&= p_2p_1'+(p_1+p_2)p_2'\,.
\end{aligned}\right.
\end{equation}

\n Set
$$\Lambda(p)\doteq\begin{pmatrix} p_1+p_2 & p_1\\\
p_2 & p_1+p_2 \end{pmatrix},
\qquad\qquad\Delta(p)\doteq\det\,\Lambda(p)\,,$$

\n and notice that

\begin{equation}\label{eq:delta_est}
{1\over 2}\,(p_1^2+p_2^2) \leq \Delta(p)\leq
2(p_1^2+p_2^2)\,.
\end{equation}

\n In particular, $\Delta(p)>0$ for all $p=(p_1,p_2)\not= (0,0)$.
Hence, $\Lambda(p)$ is invertible outside the origin and, for
$p\neq(0,0)$, we can restrict the study to the equivalent system

\begin{equation}
\left\{\begin{aligned} p'_1 &=\Delta(p)^{-1}\big[- p_1^2
+(h'_1-h'_2)
p_1+h'_1p_2\big]\,,\\
p'_2 &=\Delta(p)^{-1}\big[-p_2^2 +(h'_2-h'_1)
p_2+h'_2p_1 \big]\,.
\end{aligned}\right.
\end{equation}

Now define a new variable $s$ such that $ds/dx= \Delta(p)^{-1}$.
Using $s$ as a new independent variable, we write $p_i=p_i(s)$ and
$h_i=h_i(x(s))$ and study the equivalent system
\begin{equation}\label{eq:HJ_expl}
\left\{\begin{aligned}
{d\over ds} p_1 &=(h'_1-h'_2) p_1+h'_1p_2-p_1^2\,,\\
&\\
{d\over ds} p_2 &=(h'_2-h'_1) p_2+h'_2p_1-p_2^2 \,.\\
\end{aligned}\right.
\end{equation}

We underline that it is possible to choose the rescaling in order to
map $0$ to $0$. This choice will be assumed in the following, so
that $s(0)=0$.

\n In this new variable, as it was proved~\cite{BressanPriuli},
every unbounded trajectory $p(s)$ of~\eqref{eq:HJ_expl} actually
blows up at finite $s_o$, and it corresponds to an unbounded
trajectory $p(x)$ that tends to $\infty$ as $|x|\to\infty$. Since

$${\left|dx\over ds\right|}= \Delta\big(p(s)\big) \geq
{c_o\over (s_o-s)^2}\,,$$

\n it follows that $u(x)$ increases more than linearly as
$x\to\infty$. Therefore, $u$ is not admissible.

\n It remains to consider trajectories of~\eqref{eq:HJ_expl} that
tend to the origin, i.e. to the point where our change of variables
is singular. In~\cite{BressanPriuli} it was proven that,
by~\eqref{eq:delta_est}, these solutions satisfy
$$
\left|{dx\over ds}\right|= \Delta\big(p(s)\big)=\Og(1)\cdot e^{-2
c_o |s|}\,.
$$
\n In the original variable $x$, to the whole trajectory $s\mapsto
p(s)$ there corresponds only a portion of trajectory $x\mapsto
p(x)$, say either for $x\in ~]x_o,\infty[$ or $x\in~]-\infty, x_o[$.
Another trajectory $s\mapsto \hat p(s)$ has to be constructed to
extend the solution to all $x\in\R$.

\n For the system \eqref{eq:HJ_expl}, in the case of smooth
functions $h_1, h_2$ such that $|h_i'(x)|\leq C$, we already know
the following results (see~\cite{BressanPriuli}):

\begin{thm}\label{thm:res_coop_old} Let the cost functions $h_1, h_2$ be
smooth, and assume that their derivatives satisfy
$${1\over C}\leq h_i'(x)\leq C$$
for some constant $C>1$ and all $x\in\R$. Then the
system~\eqref{eq:HJ_impl} has an admissible solution and the
corresponding functions $\alpha_i^*=-u'_i$ provide a Nash
equilibrium solution to the non-cooperative
game~\eqref{eq:dynamics}-\eqref{eq:init_state}. Assume moreover that
the oscillation of their derivatives satisfies
$$
\sup_{x,y\in\R} \big|h_i'(x)-h_i'(y)\big|\leq\delta\,, \qquad i=1,2
$$
for some $\delta>0$ sufficiently small (depending only on $C$). Then
the admissible solution is also unique.
\end{thm}

\begin{thm}\label{thm:res_conf_old} Let any two constants
$\kappa_1,\kappa_2$ be given, with
$$\kappa_1<0<\kappa_2\,,\qquad\qquad \kappa_1+\kappa_2\not= 0\,.
$$ Then there exists $\delta>0$ such that the following
holds. If  $\,h_1,h_2$ are smooth functions whose derivatives
satisfy
$$\big|h'_1(x)-\kappa_1\big|\leq\delta\,,\qquad\qquad
\big|h'_2(x)-\kappa_2\big|\leq\delta\,,$$ for all $x\in\R$, then the
system of H-J equations~\eqref{eq:HJ_expl} has a unique admissible
solution. \end{thm}

In this paper, we want to look for admissible solutions when
smoothness of functions $h_i$ is relaxed. Namely we consider
functions $h_i$ that are piecewise linear, with a finite number of
discontinuity in their derivatives. In other words we require that
there exists a finite subdivision
$$x_o=-\infty<x_1<\ldots<x_N<x_{N+1}=+\infty$$
of $[-\infty,+\infty]$ and two $(N+1)$-tuple of constants
$(\kappa_i^1,\ldots,\kappa_i^{N+1})$, $i=1,2$, such that

\begin{equation}\label{eq:cond_on_h}
h'_i(x)=\kappa_i^j \,\,\,\, \mbox{if }x\in~]x_j,x_{j+1}[
\qquad\qquad i=1,2,\qquad j=0,\ldots,N\,.
\end{equation}

Could be of use to remark that this assumption on $h'_i$ means that
the system~\eqref{eq:HJ_expl} follows different dynamics in each
interval $\I_j\doteq ~]x_j,x_{j+1}[$: indeed, in each
$\I_j$,~\eqref{eq:HJ_expl} will have an equilibrium in $(0,0)$ and a
second one in the point $K^j=(\kappa_1^j,\kappa_2^j)$.

We also introduce the following notation (see Figure~1)
\begin{figure}\begin{center}
\psfrag{x}{$p_1$} \psfrag{y}{$p_2$} \psfrag{4}{${\cal
A}_1$}\psfrag{1}{${\cal A}_2$}\psfrag{2}{${\cal
A}_3$}\psfrag{3}{${\cal A}_4$} \psfrag{8}{${\cal
A}_5$}\psfrag{7}{${\cal A}_6$}\psfrag{6}{${\cal
A}_7$}\psfrag{5}{${\cal A}_8$}
\includegraphics[width=8cm]{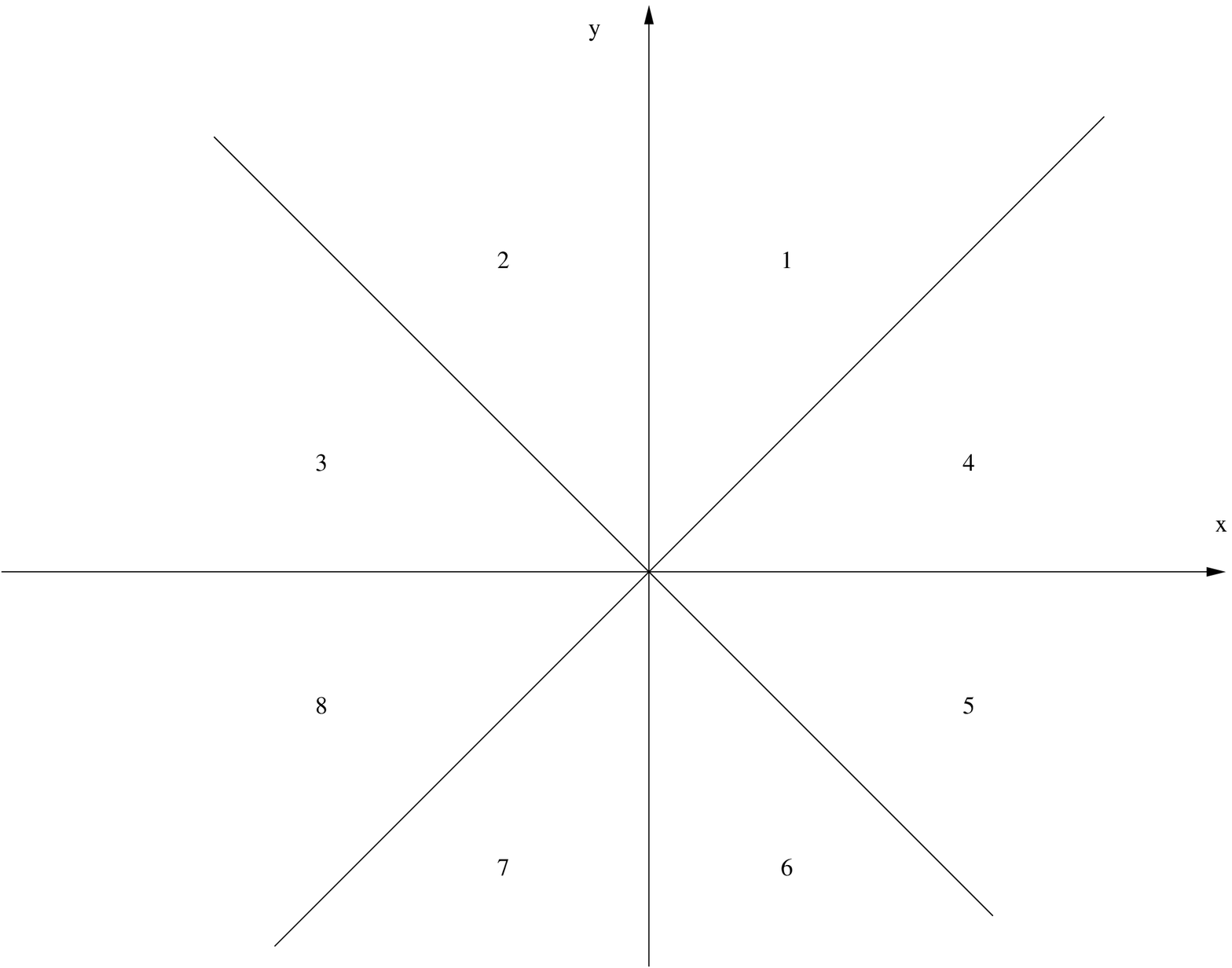}\\ Figure~1
\end{center}\end{figure}
\begin{equation}\label{eq:sets}
\mathcal{A}_i=\left\{\rho(\cos\theta,\sin\theta)\in\R^2\,\bigg|\,
\rho>0,\,\,\theta\in~\left](i-1)\disp\frac{\pi}{4},i\disp\frac{\pi}{4}\right[~\right\}\,,
\end{equation}

\n to label regions in $\R^2$, where we put our non-zero equilibria
$K^j=(\kappa_1^j,\kappa_2^j)$.

Finally, we state a couple of easy properties we will need in the
following. They provide expressions for both eigenvalues and
eigenvectors of the system obtained linearizing~\eqref{eq:HJ_expl}
around the origin. These expressions were already found
in~\cite{BressanPriuli}, and they follow from simple linear algebra.

\begin{prop}
The linearized system near $(0,0)$, corresponding
to~\eqref{eq:HJ_expl}, has the following form
\begin{equation}\label{eq:linearized}
\left(\begin{array}{c}
 p_1' \\ p_2'\\
\end{array}\right) = H \cdot
\left(\begin{array}{c}
p_1 \\ p_2\\
\end{array}\right)
\,\,, \qquad\qquad H=
\left(\begin{array}{cc}
 {\kappa_1-\kappa_2} & \kappa_1 \\ \kappa_2 & {\kappa_2-\kappa_1}\\
\end{array}\right)\,.
\end{equation}

\n Moreover the eigenvalues of the matrix $H$ are
\begin{equation}\label{eq:eigenvalues}
\lambda_-=-
\sqrt{(\kappa_1)^2+(\kappa_2)^2-\kappa_1\kappa_2}\,,\qquad
\lambda_+=
\sqrt{(\kappa_1)^2+(\kappa_2)^2-\kappa_1\kappa_2}\,,
\end{equation}

\n with corresponding eigenvectors

\begin{equation}\label{eq:eigenvectors}
\begin{array}{l}
v_- =\left(
~1, ~\displaystyle\frac{\kappa_2-\kappa_1-\sqrt{(\kappa_1)^2+(\kappa_2)^2-\kappa_1\kappa_2}}{\kappa_1}
~\right)\,,\\ \\
v_+ =\left(
~1, ~\displaystyle\frac{\kappa_2-\kappa_1+\sqrt{(\kappa_1)^2+(\kappa_2)^2-\kappa_1\kappa_2}}{\kappa_1}
~\right)\,.
\end{array}
\end{equation}
\end{prop}

\n One can immediately see that the eigenvectors
in~\eqref{eq:eigenvectors} depend actually by the ratio between
$\kappa_2$ and $\kappa_1$ only. Moreover it turns out that this kind
of dependence is indeed monotone increasing, as proved in the
following Proposition.

\begin{prop}\label{prop:alphas}
Set $\alpha=\frac{\kappa_2}{\kappa_1}$. Then the directions
corresponding to the eigenvectors $v_-$ and $v_+$ are given
(respectively) by the maps

$$
G_-(\alpha)\colon \left\{\begin{array}{rcl}
~]\,0,\infty\,[ &\to& \,]-2,-\displaystyle\frac{1}{2}\,[ \\ \\
\alpha & \mapsto & \alpha-1 -\sqrt{\alpha^2-\alpha+1}
\end{array}\right.
$$$$
g_-(\alpha)\colon \left\{\begin{array}{rcl}
~]-\infty,0\,[ & \to & \,]-\infty,-2\,[ \\ \\
\alpha & \mapsto & \alpha-1 -\sqrt{\alpha^2-\alpha+1}
\end{array}\right.
$$$$
G_+(\alpha) \colon \left\{\begin{array}{rcl}
~]\,0,\infty\,[ & \to & \,]\,0,\infty\,[ \\ \\
\alpha & \mapsto & \alpha-1 +\sqrt{\alpha^2-\alpha+1}
\end{array}\right.
$$$$
g_+(\alpha) \colon \left\{\begin{array}{rcl}
~]-\infty,0\,[ & \to & \,]-\displaystyle\frac{1}{2},0\,[ \\ \\
\alpha & \mapsto & \alpha-1 +\sqrt{\alpha^2-\alpha+1}
\end{array}\right.
$$

\n depending on the sign of $\alpha$ (and hence of
$\kappa_1\cdot\kappa_2$). These maps satisfy

\begin{equation}
\begin{array}{cc}
\frac{d}{d\alpha}G_->0\,, \qquad & \qquad \frac{d}{d\alpha}G_+>0\,,\\ \\
\frac{d}{d\alpha}g_->0\,, \qquad & \qquad \frac{d}{d\alpha}g_+>0\,.
\end{array}
\end{equation}\end{prop}

\proof The properties follow from
$$
G_-'(\alpha)=g_-'(\alpha)
= 1 -
\disp\frac{2\alpha-1}{2\sqrt{\alpha^2-\alpha+1}}=
\disp\frac{\sqrt{(2\alpha-1)^2+3}-(2\alpha-1)}{2\sqrt{\alpha^2-\alpha+1}}>0\,,
$$$$
G_+'(\alpha)=g_+'(\alpha)
= 1 + \disp\frac{2\alpha-1}{2\sqrt{\alpha^2-\alpha+1}}=
\disp\frac{\sqrt{(2\alpha-1)^2+3}+(2\alpha-1)}{2\sqrt{\alpha^2-\alpha+1}}>0\,,
$$
\n and from
$$\lim_{\alpha\to 0^+} G_-(\alpha)=\lim_{\alpha\to 0^-} g_-(\alpha)=-1-1=-2\,,$$
$$\lim_{\alpha\to 0^+} G_+(\alpha)=\lim_{\alpha\to 0^-} g_+(\alpha)=-1+1=0\,,$$
$$\lim_{\alpha\to +\infty} G_-(\alpha)=\lim_{\alpha\to +\infty} -\disp\frac{\alpha}{\alpha-1 +\sqrt{(\alpha-1)^2+\alpha}}=-\disp\frac{1}{2}\,,$$
$$\lim_{\alpha\to -\infty} g_+(\alpha)=\lim_{\alpha\to -\infty} -\disp\frac{\alpha}{\alpha-1 -\sqrt{(\alpha-1)^2+\alpha}}=-\disp\frac{1}{2}\,,$$
$$\lim_{\alpha\to +\infty} G_+(\alpha)=+\infty\,,$$
$$\lim_{\alpha\to -\infty} g_-(\alpha)=-\infty\,.$$
\qed

\section{Cooperative Situation}

We start considering all $K^j=(\kappa_1^j,\kappa_2^j)$ in ${\cal
A}_1\cup{\cal A}_2$. Notice that a similar analysis, with
straightforward adaptations, can be done if the $K^j$ are in ${\cal
A}_5\cup{\cal A}_6$. This choice implies that our system follows the
dynamics depicted in Figure~2.

\begin{figure}\begin{center}
\psfrag{x}{$p_1$} \psfrag{y}{$p_2$} \psfrag{a}{$K$}
\includegraphics[width=9cm]{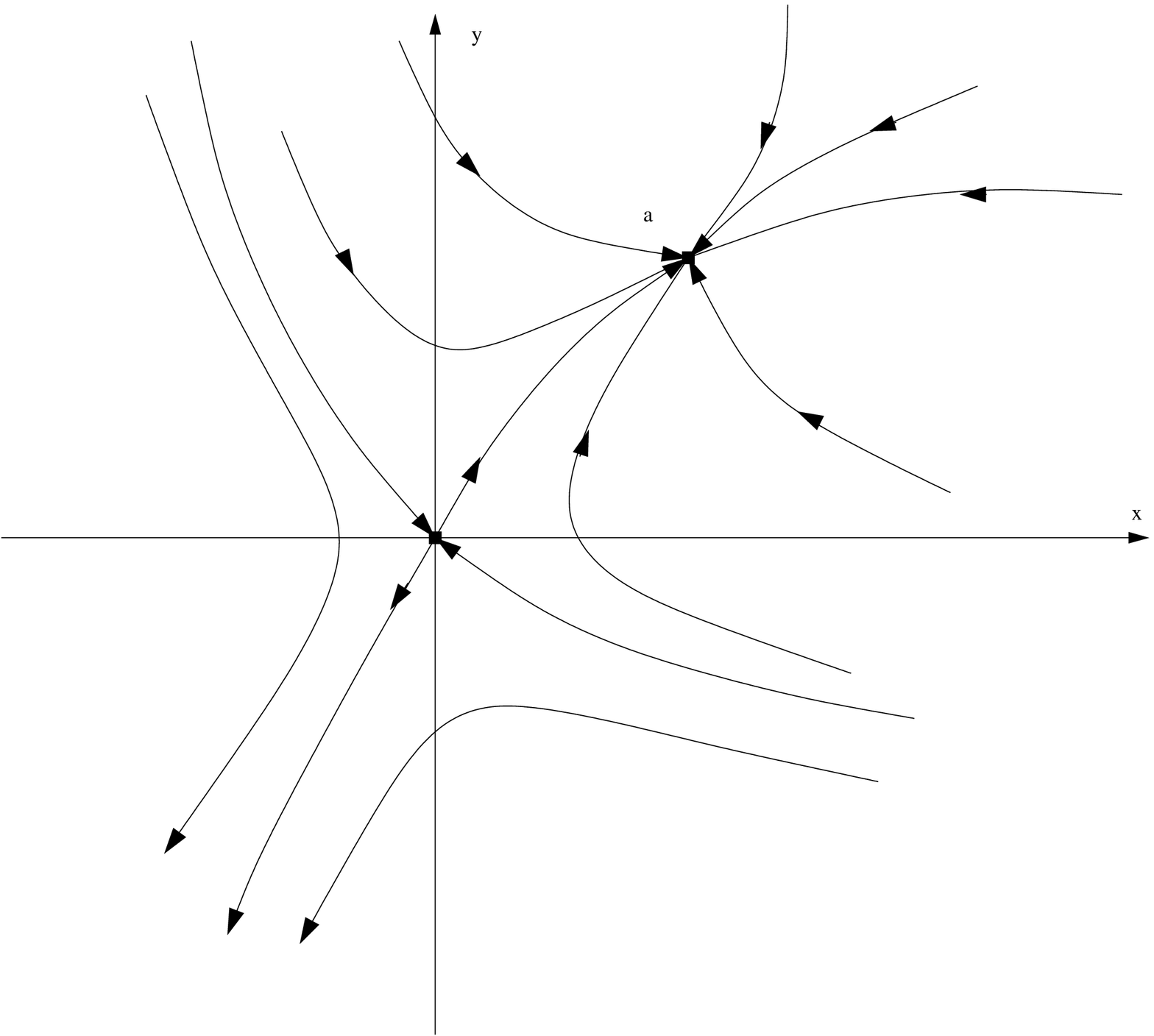} \\Figure~2
\end{center}\end{figure}

\begin{thm}\label{thm:res_coop1} Let the cost functions $h_1, h_2$ be as in~\eqref{eq:costs}, and
assume that the constants $(\kappa_1^j,\kappa_2^j)$ are all chosen
in ${\cal A}_1\cup{\cal A}_2$. Then the system~\eqref{eq:HJ_impl}
has a unique admissible solution and the corresponding functions
$\alpha_i^*=-u'_i$ provide a Nash equilibrium solution to the
non-cooperative game~\eqref{eq:dynamics}-\eqref{eq:init_state}.
\end{thm}

\proof {\bf Existence.} The existence of an admissible solution is
very easy to prove. Indeed, it is enough to glue together pieces of
admissible solutions in each interval $\I_j$. We proceed as follows:

\begin{itemize}
\item in $\I_o$, we set $p^o\equiv K^o=(\kappa_1^o,\kappa_2^o)$;
\item for $j\geq 1$, in $\I_j$ we set $p^j$ the unique solution of
the Cauchy problem for~\eqref{eq:HJ_expl} with initial datum
$p(s(x_j))=p^{j-1}(s(x_j))$. Since the set
$$
\Gamma^j=\left\{(p_1,p_2)\,\big|\, p_1,p_2\in[0,2C_1],
p_1+p_2\geq\frac{C_2}{2}\right\}\,,
$$
\n where
$$
C_1=\max\{\kappa_1^j,\kappa_2^j, \frac{p_1^{j-1}(s(x_j))}{2},
\frac{p_2^{j-1}(s(x_j))}{2} \}\,,
$$$$
C_2=\min\{\kappa_1^j,\kappa_2^j,
p_1^{j-1}(s(x_j))+p_2^{j-1}(s(x_j))\}\,,
$$
\n is
positively invariant for~\eqref{eq:HJ_expl}, each $p^j$ will exists
up to $s(x^{j+1})$ without reaching $(0,0)$ and remaining bounded;
\end{itemize}
\n Then, it is well defined the continuous function $\bar p$ given
by $\bar p(x)=p^j(x)$ whenever $x\in\I_j$. Its admissibility is an
immediate consequence of its continuity and the admissibility of
each $p^j$.\\

\n {\bf Uniqueness.} To prove that the solution built above is the
unique admissible solution to~\eqref{eq:HJ_expl}, we start proving
uniqueness on $\I_o$.

We know from~\cite{BressanPriuli} that, for $s$ negative small
enough (eventually for $s\to-\infty$), the only solutions that
remain bounded are the equilibrium $K^o$ itself and the unstable
orbits exiting from the origin. Therefore, these are the unique
possible choices, in order to retaain admissibility. If we choose an
unstable orbit in place of $K^o$, in the original variable $x$ it
would correspond to a solution defined only for $x>x_o$ (for a
suitable $x_o$). To define the solution also for $x<x_o$, we should
need a solution to
$$
\left\{\begin{array}{l}
p_1'=(\kappa^o_1-\kappa_2^o) p_1+\kappa^o_1p_2-p_1^2\,,\\
\\
p_2'=(\kappa^o_2-\kappa^o_1) p_2+\kappa^o_2p_1-p_2^2\,,
\end{array}\right.
$$
that tends to the origin as $s\to+\infty$ and remains bounded for
all negative $s$. But we know from~\cite{BressanPriuli} that no
solution with both these properties exists. Hence the uniqueness of
the solution follows on $\I_o$.

For $s>s(x_1)$, the smoothness of the right hand side
of~\eqref{eq:HJ_expl} in each interval $\I_j$ ensures that $\bar p$
is the unique continuous solution.

It remains to prove that there exists no solution with admissible
jumps in $s>s(x_1)$. But this property follows
from~\eqref{eq:jumps_2p} and from the positive invariance of the
sets
$$
\Gamma^+\doteq\left\{ (p_1,p_2)\in\R^2 \big| p_1\geq 0,\,\,p_2\geq 0
\right\}\,,
$$
$$
\Gamma^-\doteq\left\{ (p_1,p_2)\in\R^2 \big| p_1\leq 0,\,\,p_2\leq 0
\right\}\,.
$$
\n Indeed, for $s>s(x_1)$ a solution can have only jumps from
$\Gamma^+$ to $\Gamma^-$. Hence, recalling~\cite{BressanPriuli},
after a first jump the solution would be forced to remain in
$\Gamma^-$ and to tend towards $\infty$. In the $x$ variable, this
would translate into a solution $u(x)$ that grows more than linearly
as $|x|\to\infty$, and this would contradict admissibility. \qed

In light of Theorem~\ref{thm:res_coop1}, on the same line
of~\cite{BressanPriuli}, it is natural to ask whether the result
still hold for perturbations of~\eqref{eq:costs} or it fails.
Actually, we can prove the following Theorem.

\begin{thm}\label{thm:res_coop2} Let the cost functions $h'_1, h'_2$
in~\eqref{eq:costs} be smooth, and assume that:

\begin{description}
\item{(1)} their derivatives satisfy
$${1\over C}\leq h'_i(x)\leq C$$
\n for some constant $C>1$ and all $x\in\R$;
\item{(2)} on $\I_o$, the following additional assumption is
satisfied
\begin{equation}\label{eq:small_osc}
\sup_{\xi,\eta\in\I_o} \big|h'_i(\xi)-h'_i(\eta)\big|\leq\delta
\qquad i=1,2.
\end{equation}
\n for some $\delta>0$ sufficiently small (depending only on $C$).
\end{description}

\n Then the system~\eqref{eq:HJ_impl} has a unique admissible
solution.
\end{thm}

\proof We can proceed as in Theorem~\ref{thm:res_coop1}, using
Theorem~\ref{thm:res_coop_old} to deal with the perturbations.
Indeed, for $s<s(x_1)$ Theorem~\ref{thm:res_coop_old} implies that
there exists a unique admissible solution, say $p^o$. Hence, an
admissible solution on the whole real line can be built as in the
previous case: for $x\in]x_j,x_{j+1}[$, $j\geq 1$, we define
$p(x)=p^j(x)$ where $p^j$ is the unique solution
to~\eqref{eq:HJ_impl} with initial datum
$p(s(x_j))=p^{j-1}(s(x_j))$. Exactly as in
Theorem~\ref{thm:res_coop1}, this function is well defined and is a
continuous admissible solution to~\eqref{eq:HJ_expl}. Since the sets
$\Gamma^+$ and $\Gamma^-$ are still positively invariant, also
uniqueness can be proved by means of the same arguments used in
Theorem~\ref{thm:res_coop1}. \qed

\begin{rem}
We underline that the presence of the small oscillations
assumption~\eqref{eq:small_osc} is uniquely motivated by the use of
Theorem~\ref{thm:res_coop_old}, which requires~\eqref{eq:small_osc}
to provide a unique admissible solution for $s<s(x_o)$.
\end{rem}

\section{Conflicting interests}

In this section we assume that the two players have conflicting
interests, i.e. their costs satisfy $h'_1(x)\cdot h'_2(x)<-C<0$ for
all $x\in\R$. For particular choices of smooth costs, this situation
can produce infinitely many Nash equilibria to the game
(see~\cite{BressanPriuli}). Nevertheless
Theorem~\ref{thm:res_conf_old} shows that, for costs which are not
exactly opposite and under suitable assumptions of small
oscillations, it is possible to recover existence and uniqueness of
Nash equilibria. This is not the case for costs as
in~\eqref{eq:costs}.

\begin{figure}\begin{center}
\psfrag{x}{$p_1$} \psfrag{y}{$p_2$}
\psfrag{a}{$\scriptstyle{\kappa_1}$}\psfrag{b}{$\scriptstyle{\kappa_2}$}
\psfrag{u}{$\scriptstyle{\gamma_U}$}\psfrag{s}{$\scriptstyle{\gamma_S}$}
\includegraphics[width=9cm]{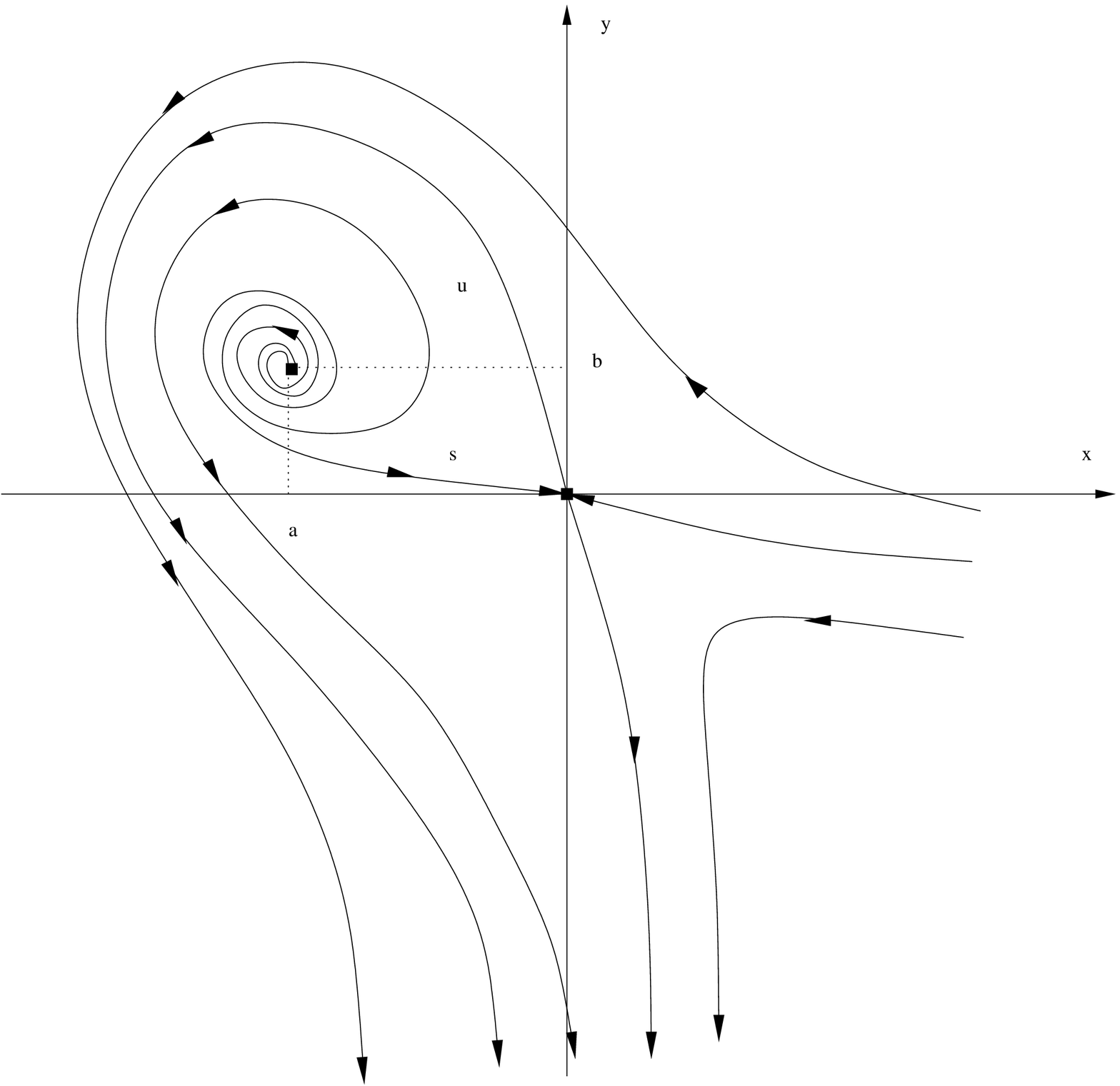} \\ Figure~3
\end{center}\end{figure}

\subsection{Case 1}

Let us consider $j=1$ in~\eqref{eq:costs}, i.e. let us consider cost
functionals that have a single jump in their derivatives. In
particular, assume this jump is located at $x=s(x)=0$. Moreover, let
us choose the constants $K^j=(\kappa_1^j,\kappa_2^j)$, $j=0,1$, so
that $K^o\in{\cal A}_4$ and $K^1\in{\cal A}_3$.

\n Under these assumptions,  the dynamics followed by the system are
depicted in Figure~3 (for $x< 0$) and Figure~4 (for $x>0$).
\begin{figure}\begin{center}
\psfrag{x}{$p_1$} \psfrag{y}{$p_2$}
\psfrag{a}{$\scriptstyle{\kappa_1}$}\psfrag{b}{$\scriptstyle{\kappa_2}$}
\psfrag{u}{$\scriptstyle{\gamma_U}$}\psfrag{s}{$\scriptstyle{\gamma_S}$}
\includegraphics[width=9cm]{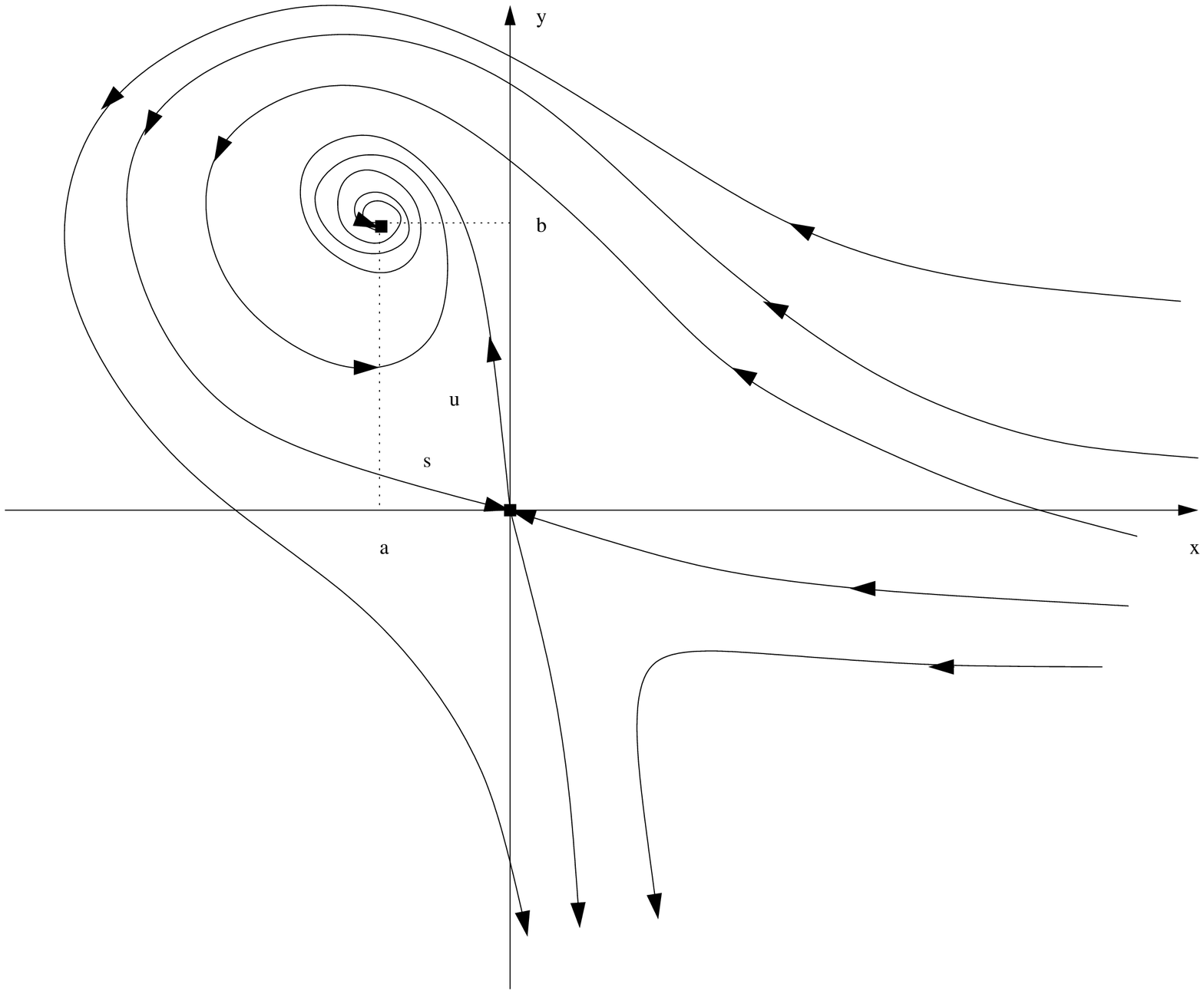} \\ Figure~4
\end{center}\end{figure}
We now prove that we could find infinitely many solutions to our
problem. Indeed, consider an initial datum
$p^\ini=(p_1^\ini,p_2^\ini)$ such that $p_1^\ini+p_2^\ini=0$ and
$p_1^\ini<0<p_2^\ini$. Recalling Proposition~\ref{prop:alphas} and
setting
$\alpha^o=\frac{\kappa_2^o}{\kappa_1^o},\alpha^1=\frac{\kappa_2^1}{\kappa_1^1}$,
we have
$$
g_-(\alpha^1)<-2<-1=\frac{p_2^\ini}{p_1^\ini}<-\frac{1}{2}<g_+(\alpha^o)\,,
$$
\n i.e. $p^\ini$ belongs to the region between the stable orbit for
the negative system (say $\gamma_S^-$) and the unstable one for the
positive system (say $\gamma_U^+$), provided it's been chosen
sufficiently near the origin. Therefore to any choice of $p^\ini$
there corresponds an admissible solution tending respectively to
either $K^1$ or $K^o$ as $s\to\pm\infty$.

\n Moreover, if the unstable orbit for the dynamics in Figure~3 (say
$\gamma_U^-$) intersects the stable one for the dynamics in Figure~4
(say $\gamma_S^+$), we can obtain an additional solution considering
as initial datum that point of intersection. Indeed the function
given by the juxtaposition of $\gamma_U^-$ and $\gamma_S^+$
corresponds, in the original variable $x$, to a solution defined on
a bounded interval $[x_-,x_+]$, with $x_-<0<x_+$ by the choice of
the rescaling. This solution can then be extended to an admissible
trajectory defined on the whole real line by using $\gamma_S^-$ for
$x<x_-$ and $\gamma_U^+$ for $x>x_+$.

\begin{rem}
The same construction can be applied when $K^o\in{\cal A}_8$ and
$K^1\in{\cal A}_7$.
\end{rem}

\subsection{Case 2}

Now we want to show, by means of a second example, how a simple
change between the positive and negative behaviors of the costs, can
lead to completely different result. Namely, we consider costs with
a single jump in their derivatives, located in $x=s(x)=0$, and
$K^o\in{\cal A}_3$, $K^1\in{\cal A}_4$. This choice produce a game
with no admissible solutions to~\eqref{eq:HJ_expl}.

\n We proceed by contradiction. Assume that an admissible solution
$\tilde p=(\tilde p_1,\tilde p_2)$ exists, for a Cauchy problem with
initial datum $\tilde p (0)=p^\ini$. Then, recalling the results
in~\cite{BressanPriuli}, we have that
$$
\lim_{s\to +\infty} |\tilde p(s)|<+\infty
$$
actually implies
$$
\lim_{s\to +\infty} |\tilde p(s)|=0\,,
$$
and hence $\tilde p$ is one of the stable orbits of the positive
system. Now we underline that this means $p^\ini\notin \gamma_U^+$.
Then, we can repeat the proof of Theorem~\ref{thm:res_conf_old},
given in~\cite{BressanPriuli}, and find
$$
\lim_{s\to s_o +} |\tilde p(s)|=+\infty\,,
$$
\n for a suitable $s_o<0$, eventually $s_o=-\infty$. Therefore the
solution cannot be admissible, and we have a contradiction.

\n Notice that the previous calculations hold even if the unstable
orbit for the dynamics in Figure~3 (say $\gamma_U^-$) intersects the
stable one for the dynamics in Figure~4 (say $\gamma_S^+$). This
means there is no solution as the one built in the previous case,
using more trajectories in the $s$ variable: this is obviously due
to the fact that we cannot find solutions bounded at $+\infty$
(resp. $-\infty$) to extend a possible $\tilde p$ when $x>x_+$
(resp. $x<x_-$).

\begin{rem}
The same result can be obtained when $K^o\in{\cal A}_7$ and
$K^1\in{\cal A}_8$.
\end{rem}

\begin{rem}
Actually, one can still construct particular cases so that there
exist admissible solutions. Fixed $K^o,K^1$ as above, assume that
the trajectories $\gamma^-_U$ and $\gamma^+_S$ intersect in a point.
Moreover, set $x_-$ and $x_+$ the values introduced in the previous
example, $\ell=|x_+-x_-|$ and ${\cal J}_n=]x_-+n\ell,x_++n\ell[$,
$n\in\Z$. We can define piecewise linear costs on the whole $\R$ by
repeating on each ${\cal J}_n$ the same $2$-value piecewise linear
cost. In other words, $\forall n\in\Z$ set

\begin{equation}
h'_i(x)_{\big| {\scriptstyle{\cal J}_n}}=
\left\{\begin{array}{ll}
\kappa_i^o &\mbox{if }x\in~]x_-+n\ell,n\ell[ \\
\kappa_i^1 &\mbox{if }x\in~]n\ell,x_++n\ell[
\end{array}\right.
\qquad\qquad i=1,2,
\end{equation}

\n Then, we find a solution by simply gluing together periodically
$\gamma^-_U$ and $\gamma^+_S$. This solution is admissible, being
bounded in the $p_1,p_2$ plane.

\n Anyway no general results as Theorem~\ref{thm:res_conf_old} is
possible.

\end{rem}

\section{Mixed Cases}

In this section we end our presentation of ill-posed problems, with
a last example presenting costs that can switch from a situation
with conflicting interests into a cooperative one. More precisely,
we consider costs with a single jump in their derivative, located
again in $x=s(x)=0$, and $K^o\in{\cal A}_5\cup{\cal A}_6$,
$K^1\in{\cal A}_1\cup{\cal A}_2$. Moreover, let us assume
\begin{equation}\label{eq:cond_ex2}
\alpha^1=\frac{\kappa_2^1}{\kappa_1^1}\neq\frac{\kappa_2^o}{\kappa_1^o}=\alpha^o\,.
\end{equation}
\n With these assumptions, the system follows the dynamics depicted
in Figure~5 (resp. Figure~2) for $x< 0$ (resp. $x>0$) and $K^o,K^1$
are not on the same line through the origin.

\begin{figure}\begin{center}
\psfrag{x}{$p_1$} \psfrag{y}{$p_2$} \psfrag{a}{$K$}
\includegraphics[width=9cm]{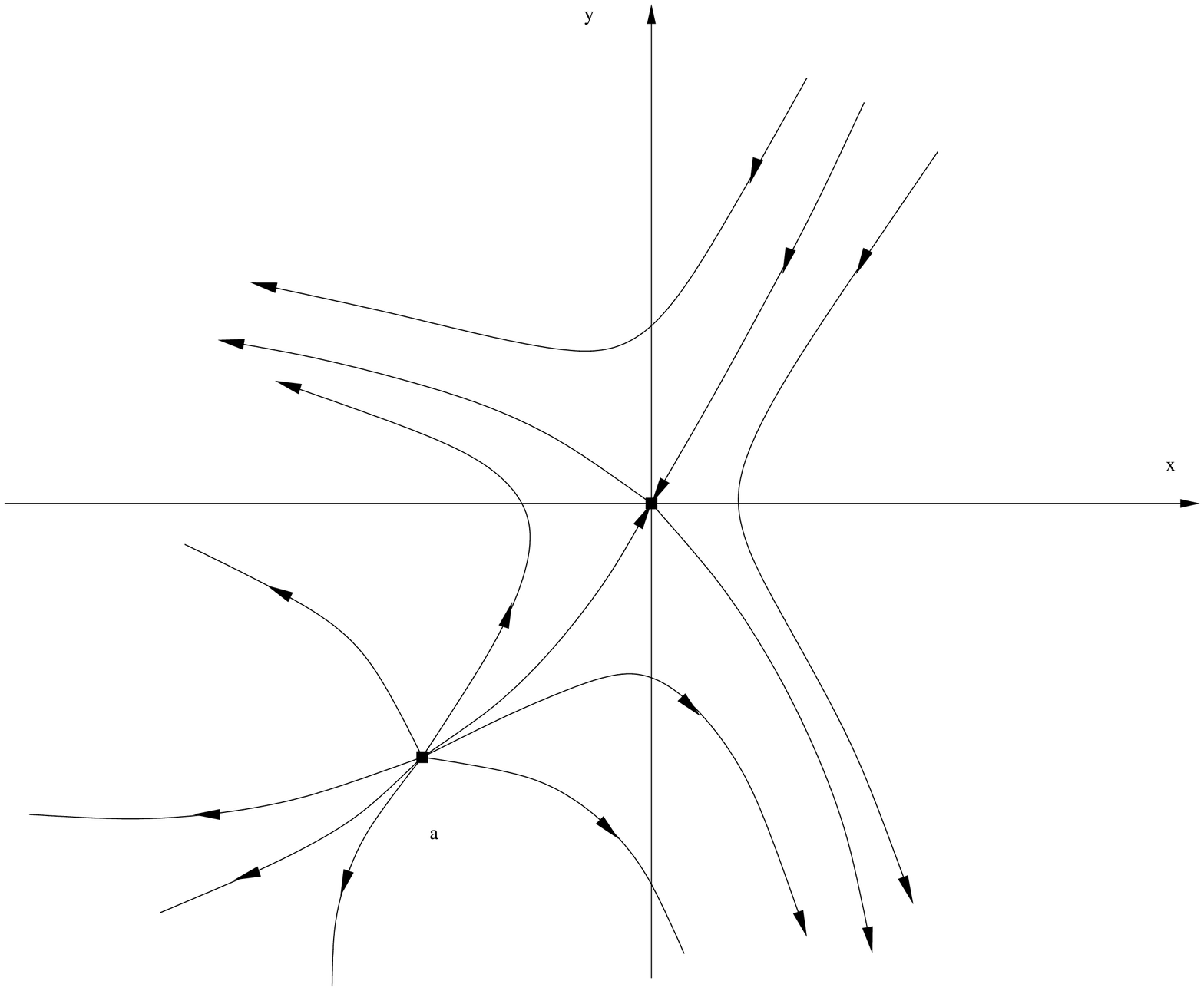} \\ Figure~5
\end{center}\end{figure}

Again, we observe the existence of infinitely many Nash equilibria.
Assume it holds $\alpha^o<\alpha^1$ in~\eqref{eq:cond_ex2} (the
opposite inequality leading to a similar analysis). Then, we can
consider the non-empty region
$$
\Omega=\left\{ (p_1,p_2)\in\R^2~ \left| ~p_1<0<p_2,~
G_-(\alpha^o)<\frac{p_2}{p_1}<G_-(\alpha^1) \right.\right\}\,.
$$
\n This region is, at least near the origin, say in a neighborhood
$\mathcal{O}$, exactly the region between the stable orbit for the
positive system and the unstable one for the negative system. Taking
as initial datum any point $p^\ini$ both in $\Omega$ and in
$\mathcal{O}$, we can construct an admissible solution in the
following way. We take for $s<0$ the unique solution to the negative
system, passing through $p^\ini$ at $s=0$ and tending to $K^o$ as
$s\to -\infty$. In an analogous way, we take for $s>0$ the unique
solution to the positive system, passing through $p^\ini$ at $s=0$
and tending to $K^1$ as $s\to +\infty$. Every such a solution, being
continuous and bounded in $s$, corresponds to an admissible solution
$u(x)$.

\end{document}